\documentclass{article}
\usepackage[T2A]{fontenc}
\usepackage[utf8]{inputenc}
\usepackage[english]{babel}
\usepackage{amsmath}
\usepackage{amsfonts}
\usepackage{amssymb}
\usepackage{wasysym}
\usepackage{amsthm}
\usepackage{graphicx}
\usepackage{setspace}
\usepackage{tikz}
\usepackage{float}
\usepackage{tikz-cd}
\usetikzlibrary{decorations.markings}
\usepackage{ytableau}
\usepackage{longtable}
\usepackage{empheq}
\bibliography{references}
\usepackage[hidelinks,colorlinks=true,unicode]{hyperref}
\hypersetup{
	linkcolor=blue,
	citecolor=blue,
	filecolor=magenta,
	urlcolor=blue,
}
\usepackage{url}
\usepackage[left=2cm,right=2cm,top=2cm,bottom=2cm,bindingoffset=0cm]{geometry}
\usepackage[colorinlistoftodos,prependcaption,textsize=tiny]{todonotes}

\newtheorem{theorem}{Theorem}
\newtheorem{st}{Proposition}

\newtheorem{definition}{Definition}
\newtheorem{statement}{Statement}
\theoremstyle{remark}

\renewcommand{\qed}

\title{Tropical Ptolemy Transformations and Invariants of Braids} \vspace{1.2cm}
\author{Gurnoor Singh}

\begin{document}
\maketitle


\begin{abstract}
It often happens in mathematics that one and the same equation is known under different names in different areas of mathematics. The famous pentagon identity appears in low-dimensional topology in different ways. In this paper, we use the tropical version of the Ptolemy equation to construct invariants of braids.
\end{abstract}



\maketitle
\section{Introduction}

In the study of invariants, one examines mathematical objects—such as \textit{manifolds}—represented in various states, such as \textit{triangulations}, and investigates their behavior under specific local transformations. To each state, one associates a set of \textit{data}, which may consist of quantities such as \textit{lengths}, \textit{areas}, or \textit{angles}. 

To construct invariants under these transformations, as explored in \cite{Invariants-and-Pictures, Manturov-Kim-2019,Manturov-Nikonov-2023}, we treat the associated data (e.g., edge lengths) as variables and impose algebraic relations (or transformation laws) that reflect how this data evolves under local moves. A fundamental requirement is that the composition of data transformations corresponding to a sequence of local moves—such as \textit{Pachner moves}—returns the system to its original state. That is, the cumulative effect of such transformations must act as the identity on the data.

A classical instance of this phenomenon is the \textit{Ptolemy transformation}: Ptolemy’s theorem is a classical result in Euclidean geometry that characterizes cyclic quadrilaterals—quadrilaterals whose vertices lie on a common circle. The theorem provides a precise algebraic relationship between the side lengths and diagonals of such a quadrilateral.

\begin{statement}
    Let \(ABCD\) be a quadrilateral inscribed in a circle (i.e., a cyclic quadrilateral). Denote the lengths of the sides \(AB = a\), \(BC = b\), \(CD = c\), \(DA = d\), and the diagonals \(AC = x\), \(BD = y\). Then the following identity holds:
    \[
    ac + bd = xy.
    \]
\end{statement}

Ptolemy's identity appears naturally in the study of triangulations, cluster algebras, and representations of the braid group. In particular, it forms the foundation of the \emph{Ptolemy transformation}, a rule that governs how edge data—such as lengths or weights—evolve under a \emph{diagonal flip} in a triangulated polygon. Such flips replace one diagonal of a convex quadrilateral with the other, inducing a transformation consistent with Ptolemy’s relation.

\begin{figure}[h]
    \centering
    \includegraphics[width=0.5\linewidth]{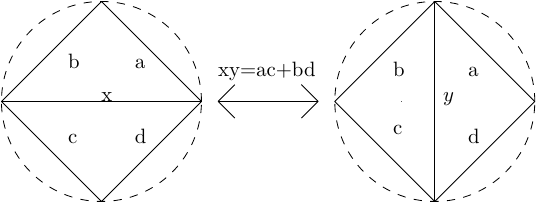}
    \caption{The Ptolemy relation in a cyclic quadrilateral: \( xy = ac + bd \)}
    \label{fig:The Ptolemy Relation}
\end{figure}

The Ptolemy transformation satisfies the celebrated \textit{pentagon identity}: the composition of five Ptolemy transformations, arranged in a specific sequence corresponding to five flips of diagonals in a pentagon, yields the identity transformation. This structure has been studied extensively in the literature; see, for example, \cite{Felikson}. Geometrically, if two quadrilaterals share three vertices and both are cyclic, then the entire pentagon is also cyclic. The data associated with any one of its quadrilaterals can then be recovered from the data of the other two.

A matrix-theoretic analogue of this idea was recently explored in \cite{Rohozhkin-2025}, where matrices were defined to encode the movement of points on the standard sphere, using Delaunay triangulations. From this construction, an invariant of the pure braid group was obtained.

Building upon the framework of invariants for pure braid groups derived from general Ptolemy transformations, this paper demonstrates that the \textit{tropical Ptolemy transformation} also satisfies the pentagon identity. We employ the tropical version of the identity, as studied in \cite{MS}. This leads to the construction of new invariants for braid groups, thereby extending the applicability of Ptolemy-type transformations into the realm of tropical geometry. Our work thus provides novel insights into braid invariants through this tropical lens.

The present paper is organized as follows: In Section 2, we recall the notion of Delaunay triangulations on the 2-sphere and flips on them. We study their essential properties, called Pentagon relation and far-commutativity. In Section 3, we define a labeling of edges of Delaunay triangulations valued in a tropical semifield and correspondence between labelings of two Delaunay triangulations related by a flip. In Section 4, a sequence of Delaunay triangulations and flips is induced by a given braid. In Section 5, an invariant for braids is constructed by using the associated sequence of Delaunay triangulations and flips, and labelings valued in Tropical semifield.

\section{Delaunay Triangulations on the sphere}
We follow \cite{Invariants-and-Pictures}.

\begin{definition}

Let \( z_i(t) \in \mathbb{S}^2 \), for \( i = 1, \ldots, n \) and \( t \in [0,1] \), be continuous trajectories of \( n \) moving points on the sphere. At each time \( t \), define the corresponding point configuration
\[
z(t) := \{ z_1(t), z_2(t), \ldots, z_n(t) \} \subset \mathbb{S}^2.
\]

We define the \emph{Voronoi cell} associated to the point \( z_i(t) \) as:
\[
U_i(t) := \left\{ z \in \mathbb{S}^2 \;\middle|\; \forall j \neq i, \; \| z - z_i(t) \| \leq \| z - z_j(t) \| \right\}.
\]

\end{definition}
The collection \( \{ U_i(t) \}_{i=1}^n \) forms a tiling of the plane by closed, convex regions. The boundaries between the regions are formed by segments of perpendicular bisectors between pairs of points, and the union of these boundaries defines a planar graph \( \Gamma(t) \), called the \emph{Voronoi diagram} associated to \( z(t) \). Generically, \( \Gamma(t) \) is a trivalent graph with some unbounded (infinite) edges.

\begin{definition}[General Position on the Sphere]
Let \( \mathbb{S}^d \) denote the standard unit \( d \)-sphere, and let \( P = \{p_1, p_2, \ldots, p_n\} \subset \mathbb{S}^d \) be a finite collection of points on the sphere.
We say that the set \( P \subset \mathbb{S}^d \) is in \emph{general position} if the following condition holds:

For every integer \( k \in \{1, \ldots, d+2\} \), no subset of \( k \) points from \( P \) lies entirely on a great \( (k-2) \)-dimensional subsphere of \( \mathbb{S}^d \).

\noindent Equivalently:
\begin{itemize}
    \item No two points coincide.
    \item No three points lie on a great circle (i.e., a 1-dimensional subsphere).
    \item No four points lie on a common great 2-sphere.
    \item More generally, for each \( k \in \{2, \ldots, d+1\} \), no \( k+1 \) points lie on a great \( k \)-sphere.
\end{itemize}
\end{definition}






\begin{definition}

The \emph{Delaunay triangulation} \( D_t \) is the planar dual of the Voronoi diagram \( \Gamma(t) \): vertices of \( D_t \) correspond to the points \( z_i(t) \), and edges are drawn between two vertices if their Voronoi cells share an edge. Generically, \( D_t \) consists of triangles.

A configuration \( z(t) \) is said to be \emph{generic} if no four points lie on a common circle and no three are collinear. As \( t \) varies, the diagram \( \Gamma(t) \) evolves continuously, except at finitely many singular moments \( t = t_k' \in [0,1] \), where a codimension-one degeneracy occurs.
(see \cite{Aurenhammer} for a comprehensive discussion of general position in Delaunay triangulations).

\end{definition}

\begin{figure}
    \centering
    \includegraphics[width=0.4\linewidth]{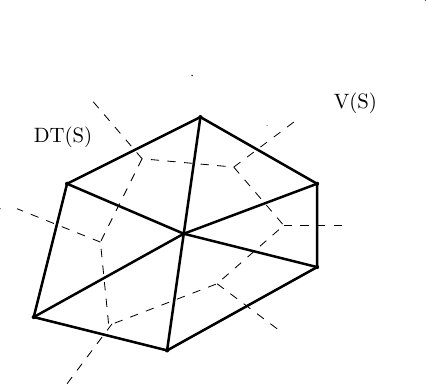}
    \caption{Voronoi Diagram and Delaunay Triangulation}
    \label{fig:Voronoi Diagram and Delaunay Triangulation}
\end{figure}

\subsection*{Euler Characteristic for Delaunay Triangulation on the Sphere}
Let \( P = \{p_1, p_2, \dots, p_n\} \subset \mathbb{S}^2 \) be a set of \( n \) points on the sphere. The Delaunay triangulation \( \mathcal{T} \) of \( P \) induces a planar graph \( G = (V, E, F) \), where:
\begin{itemize}
    \item \( V = n \) is the number of vertices (points in \( P \)),
    \item \( E \) is the number of edges,
    \item \( F \) is the number of faces (triangles in the triangulation).
\end{itemize}

The Euler characteristic of the sphere \( \mathbb{S}^2 \) is \( \chi(\mathbb{S}^2) = 2 \).
 The relationship between edges and faces in a triangulation is \( E = \frac{3F}{2} \), leading to:
\[
F = 2n - 4 \quad \text{and} \quad E = 3n - 6.
\]
Thus, the number of edges in the Delaunay triangulation is \( E = 3n - 6 \).
\begin{figure}
    \centering
    \includegraphics[width=0.5\linewidth]{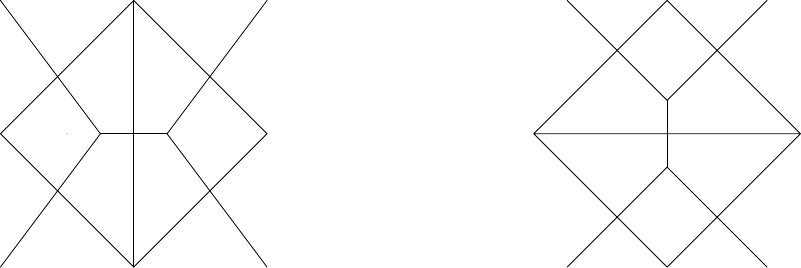}
    \caption{Voronoi Tiling Change}
    \label{fig:Voronoi Tiling Change}
\end{figure}
\subsection*{Flips}

For a finite set \( P \subset \mathbb{S}^2\), if four points lie on a common circle, then Delaunay triangulation is not determined uniquely.
\begin{definition}
Let \( T \) be a triangulation of a convex polygon \( P \subset \mathbb{S}^2 \), and let \( d \in T \) be a diagonal contained in exactly two adjacent triangles \( \Delta_1, \Delta_2 \subset T \), such that their union \( \Delta_1 \cup \Delta_2 \) forms a convex quadrilateral \( Q \subset P \). Let \( \{v_i, v_j, v_k, v_\ell\} \subset V \) denote the four vertices of \( Q \) in cyclic order, with \( d \) being one diagonal of \( Q \). Then there exists a unique other diagonal \( d' \neq d \) that lies in the interior of \( Q \), connecting the other pair of non-adjacent vertices.

The operation that replaces \( d \) with \( d' \), producing a new triangulation
\[
T' := (T \setminus \{d\}) \cup \{d'\},
\]
is called a \emph{\textit{flip}}. We denote this as \( T \xrightarrow{d \mapsto d'} T' \).
\end{definition}

\begin{figure}
    \centering
\includegraphics[width=0.5\linewidth]{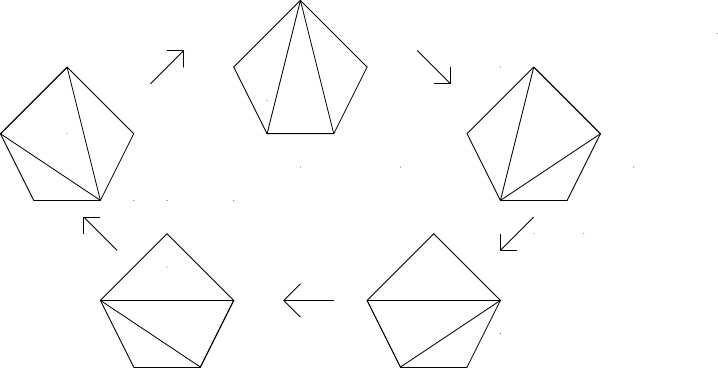}
    \caption{Pentagon Relation}
    \label{fig:Pentagon Relation}
\end{figure}
Let $P(t) = \{p_{i}(t)\}_{i=1}^{N} \subset \mathcal{S}^2$ such that $p_{i} : [0,1] \rightarrow \mathcal{S}^2$ with $p_{i}(0) = p_{i}(1)$ and $p_{i}(t) \neq p_{j}(t)$ for $i\neq j$ and for any $t \in [0,1]$. The flip operation happens in the moment $s \in [0,1]$ when $P(s)$ is not in general position, that is, four points \( \{p_i(s), p_j(s), p_k(s), p_\ell(s)\} \subset \mathbb{S}^2 \) lie on a common circle, i.e., they are \emph{co-circular} in the spherical sense. In this case, the spherical quadrilateral \( Q \subset \mathbb{S}^2 \), formed by the union \( \Delta_1 \cup \Delta_2 \), is cyclic and admits exactly two distinct diagonals $d$ and $d'$ that lie entirely within the interior of \( Q \). 
The replacement of one diagonal \( d \in T \) with the other diagonal \( d' \) results in a new collection of geodesic arcs \( T' := (T \setminus \{d\}) \cup \{d'\} \), which again forms of a triangulation of \( P \): maximality, non-crossing, and coverage of the interior of \( P \) by spherical triangles with disjoint interiors.

Moreover, the flip corresponds to a local change in the triangulation induced by a codimension-one degeneracy in the configuration space of vertex placements on \( \mathbb{S}^2 \). Specifically, such a degeneracy occurs precisely when a continuous motion of the vertices causes four of them to become co-circular (i.e., lie on a common great circle). At this critical configuration, the spherical Delaunay triangulation changes combinatorially by replacing the diagonal \( d \) of the cyclic quadrilateral \( Q \) with the unique other diagonal \( d' \), yielding a new Delaunay triangulation.

\subsection*{Pentagon Relation and Flips}

Let \( P \) be a convex pentagon with vertices \( a, b, c, d, e \) labeled cyclically. A triangulation of \( P \) consists of three non-crossing diagonals that divide the interior into three triangles.

Let \( f_i: T_i \to T_{i+1} \) denote the flip transforming triangulation \( T_i \) into \( T_{i+1} \). Then the composition of the five flips returns to the original triangulation: There are exactly five distinct triangulations \( T_0, T_1, T_2, T_3, T_4 \) of the convex pentagon, such that each \( T_{i+1} \) is obtained from \( T_i \) by a single flip (indices modulo 5). Specifically, applying the sequence of flips \( f_{0}, f_{1}, f_{2}, f_{3}, f_{4} \) in order results in:
\[
f_{4} \circ f_{3} \circ f_{2} \circ f_{1} \circ f_{0} = \mathrm{id},
\]
where \( \mathrm{id} \) denotes the identity transformation on the space of triangulations of the pentagon.

This equality is known as the \emph{pentagon relation}. It reflects the fact that performing the sequence of five flips, each modifying a single diagonal in the triangulation, results in the initial triangulation. The space of triangulations of a convex pentagon, together with flips as edges, forms the 1-skeleton of a pentagon (as a graph), and this cyclic structure gives rise to the relation above.

\subsection*{Far-Commutativity of Flips}
Let \( P \) be a convex polygon, and let \( T \) be a triangulation of \( P \). Suppose that \( Q_1 \) and \( Q_2 \) are two convex quadrilaterals contained in \( P \), each formed by the union of two adjacent triangles in \( T \), and that the interiors of \( Q_1 \) and \( Q_2 \) are disjoint. Denote by \( d_{1} \) and \( d_{2} \) the diagonals of \( Q_1 \) and \( Q_2 \) that are present in the triangulation \( T \), and let \( d_1', d_2' \) be the respective diagonals that replace them under a flip.

Define the flip operations:
\[
f_1(T) := (T \setminus \{d_1\}) \cup \{d_1'\}, \qquad f_2(T) := (T \setminus \{d_2\}) \cup \{d_2'\}.
\]

Because the interiors of \( Q_1 \) and \( Q_2 \) are disjoint, the flip \( f_1 \) only modifies edges and triangles within \( Q_1 \), and \( f_2 \) only modifies those within \( Q_2 \). Therefore, the operations \( f_1 \) and \( f_2 \) affect disjoint parts of the triangulation, and hence they commute:
\[
f_1(f_2(T)) = f_2(f_1(T)).
\]

\section{Tropical semifield  and Labels of Edges}
      
\begin{definition}
    
The \textbf{tropical semifield} is an ordered semifield equipped with two operations, \( \oplus \) and \( \otimes \), defined as:
\[
x \oplus y = \max(x, y), \quad x \otimes y = x + y.
\]
We call \( \oplus \) the tropical addition, and \( \otimes \) the tropical multiplication.

In this framework, the classical Ptolemy relation \( xy = ab + cd \) is replaced by its tropical analogue:
\begin{equation}\label{eqn:tropical_ptolemy}
x\oplus y= (a\oplus c)\otimes(b\oplus d).
\end{equation}
Equation (\ref{eqn:tropical_ptolemy}) is called the {\em tropical Ptolemy relation.}
\end{definition}
See \cite{MS} for a detailed exposition on Tropical Geometry.

\begin{figure}
    \centering
    \includegraphics[width=0.6\linewidth]{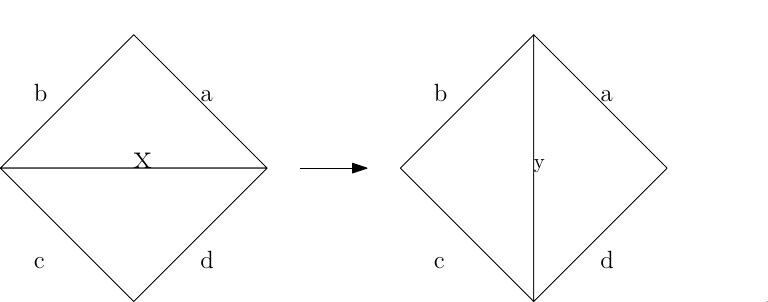}
    \caption{Tropical Ptolemy Transformation $x\oplus y= (a\oplus c)\otimes(b\oplus d)$ }
    \label{fig: Flip}
\end{figure}

\begin{definition}

    Let $T$ be a triangulation. A {\bf labeling of edges of $T$} is an assignment of elements of a tropical semifield $(X,\oplus, \otimes)$ to edges of $T$. Under a flip $f$ on a quadrilateral with labels $a,b,c,d$ and a diagonal with a label $x$, the diagonal with a label $x$ is replaced by another diagonal with a label $y$ satisfying the tropical Ptolemy relation (\ref{eqn:tropical_ptolemy}), see Fig. \ref{fig: Flip}. We call the map $F$ from labels of edges of $T$ to labels of edges of $f(T)$ such that $F(x) =y$, but labels of other edges remain, a {\bf flip of labels of edges of $T$} (or simply {\bf a flip on labels}).
\end{definition}

\subsection{Consistency with Involution, Far-Commutativity, and Pentagon Relation}\label{subsec:consistency-labels}

\begin{figure}
    \centering\includegraphics[width=0.9\linewidth]{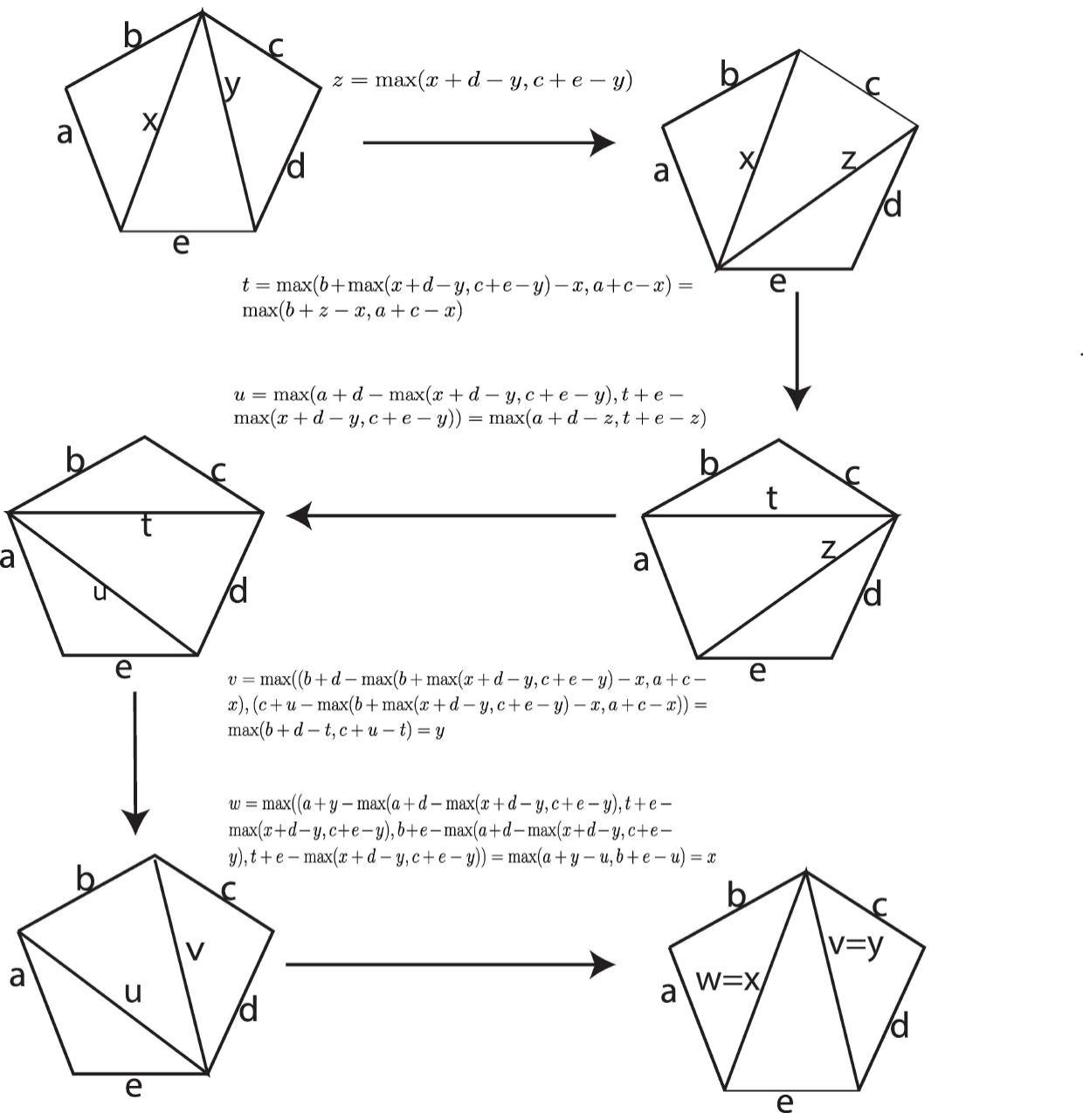}
    \caption{The pentagon equation is satisfied by tropical Ptolemy transformation. }
    \label{fig:The pentagon equation is satsfied by tropical transformation. }
\end{figure}

Now, we analyze how the flip on label is compatible with three key properties of flips of triangulations: involution, far-commutativity, and the pentagon relation. 
\begin{itemize}
    \item 

Let $f$ be a flip on a triangulation $T$ replacing a diagonal $e$ by $e'$ of a convex quadrilateral with boundary edges. Assume that edges of $T$ are labeled by elements in a tropical semifield $(X, \oplus, \otimes)$. Suppose that the boundary edges labeled \( a, b, c, d \in X \), and diagonal \( e \) labeled by \( x \in X \). Let $F$ be the flip of labels of $T$ such that  $F(x) = x'$ where $ x'= \max(a + c,\; b + d) - x$.

Applying the flip operation once more to the same quadrilateral (now with diagonal \( e' \) and same edge labels \( a, b, c, d \)) reverts the triangulation to its original state. Let \( x'' \in T \) denote the label of the diagonal \( e \) after this second flip. Using the tropical Ptolemy relation again:
\begin{equation}
x'' = \max(a + c,\; b + d) - x'. \tag{2}
\end{equation}

Substitute the expression for \( x' \) from equation (1) into equation (2):
\begin{align*}
x'' &= \max(a + c,\; b + d) - \left( \max(a + c,\; b + d) - x \right) \\
    &= \left( \max(a + c,\; b + d) - \max(a + c,\; b + d) \right) + x \\
    &= 0 + x \\
    &= x.
\end{align*}
Hence we obtain that 
\[
F(F(x)) = x.
\]
It follows that the flip \( F \) on labels is an involution on the space of tropical edge labelings:
\[
F^2 = \mathrm{id}.
\]
\item
For far-commutativity, let us consider two disjoint quadrilaterals \( Q_1 \) and \( Q_2 \) embedded in a triangulation of a surface (e.g., the sphere \( S^2 \)), such that:

\begin{itemize}
    \item \( Q_1 \) has boundary edges labeled \( a, b, c, d \in T \) and diagonal \( x \in T \),
    \item \( Q_2 \) has boundary edges labeled \( a', b', c', d' \in T \) and diagonal \( x' \in T \),
    \item \( Q_1 \cap Q_2 = \emptyset \); that is, the quadrilaterals are disjoint.
\end{itemize}

Let \( F_i \) denote the flip on labels corresponding to flip $f_{i}$ in \( Q_i \) for $i=1,2$. Then we obtain that

\begin{align}
F_1(x) &= \max(a + c, b + d) - x,\\
F_2(x') &= \max(a' + c', b' + d') - x'.
\end{align}

Since $Q_{1}$ and $Q_{2}$ are disjoint quadrilaterals, one can easily see that 
\begin{align}
F_2(F_1(x)) = F_1(F_2(x)) &= \max(a + c, b + d) - x,\\
F_2(F_1(x')) = F_1(F_2(x')) &= \max(a' + c', b' + d') - x',
\end{align}
and hence $F_{1}\circ F_{2}= F_{2}\circ F_{1}$.


\item
Now let us consider the pentagon relation. Let the edge labels involved in a sequence of flips be denoted by elements \( x, y, z, t, u, v, w \in T \). The transformation of the labels under flips is governed by the tropical Ptolemy relations:

\begin{align}
z &= \max(x + d - y,\; c + e - y), \tag{2} \\
t &= \max(b + z - x,\; a + e - x), \tag{3} \\
u &= \max(a + d - z,\; t + e - z), \tag{4} \\
v &= \max(b + d - t,\; c + u - t), \tag{5} \\
w &= \max(a + y - u,\; e + b - u). \tag{6}
\end{align}

Each equation updates the label of a diagonal or edge as a tropical maximum over sums of neighboring edge labels from the quadrilateral undergoing a flip.

Let \( T_0, T_1, T_2, T_3, T_4 \) be the five triangulations of the convex pentagon obtained by successive flips, and let \( f_i: T_i \to T_{i+1} \) denote the tropical label update at step \( i \), governed by relations (2)--(6).

Let \( A_0 \in T^N \) denote the vector of tropical edge labels in the triangulation \( T_0 \), where \( N = 3n - 6 \) is the number of edges in a triangulation of \( n \) points. Then, define:
\[
A_{i+1} = f_i(A_i), \quad \text{for } i = 0, 1, 2, 3, 4 \mod 5.
\]

The tropical pentagon relation asserts:
\[
f_4 \circ f_3 \circ f_2 \circ f_1 \circ f_0 (A_0) = A_0,
\]
see Fig. \ref{fig:The pentagon equation is satsfied by tropical transformation. } In other words, after performing the sequence of five flips on a cocircular configuration of five points, the edge labels return to their initial values (modulo relabeling consistent with tropical Ptolemy transformations).
\end{itemize}

\section{A sequence of Delaunay triangulations and flips induced by a braid}

\begin{definition}

    The \emph{spherical braid group} on \( n \) strands, denoted \( B_n(\mathbb{S}^2) \), is defined as the fundamental group of the unordered configuration space $\mathrm{UConf}_n(\mathbb{S}^2)$ of \( n \) distinct points on the 2-sphere:
\[
B_n(\mathbb{S}^2) = \pi_1(\mathrm{UConf}_n(\mathbb{S}^2)).
\]

\end{definition}

   This represents isotopy classes of $n$ points moving on the sphere $S^2$ without collisions, where the points are indistinguishable.
   For details on braid groups and their geometric representations, see \cite{Manturov-Kim-2019, Manturov-Nikonov-2023, Fedoseev-Manturov-Nikonov-2021}.

Let \( b \in B_n(\mathbb{S}^{2}) \) be represented by a continuous loop
\[
(x_1(t), \dots, x_n(t)) \in \mathrm{UConf}_n(\mathbb{S}^2), \quad t \in [0,1],
\]
where \( x_i(0) = x_i(1) \) for all \( 1 \leq i \leq n \). We further assume that the path \( (x_1(t), \dots, x_n(t)) \) is \emph{generic}, in the following precise sense:

\begin{definition}
A path \( (x_1(t), \dots, x_n(t)) \in U\mathrm{Conf}_n(\mathbb{S}^2) \), \( t \in [0,1] \), is \emph{generic} if the configuration \( \{x_1(t), \dots, x_n(t)\}\subset\mathbb{S}^2\) satisfies:
\begin{itemize}
    \item For all but finitely many \( t \in [0,1] \), the points are in general position—i.e., no four points lie on a common circle.
    \item At finitely many critical times \( 0 < t_1 < t_2 < \cdots < t_\ell < 1 \), exactly four of the points lie on a circle $C$ such that there is no points in one of connected components of $\mathbb{S}^{2}\backslash C$.
\end{itemize}
\end{definition}

For a generic braid $b= \{x_1(t), \dots, x_n(t)\}$, one can obtain Delaunay triangulation $T(t)$ of $\mathbb{S}^{2}$ with $\{x_1(t), \dots, x_n(t)\}$ for $t \in [0,1]$ except finitely many times \( 0 < t_1 < t_2 < \cdots < t_\ell < 1 \). By definition of generic braid, $T(s_{1})$ and $T(s_{2})$ for $s_{1}, s_{2} \in (t_{i},t_{i+1})$ for $i=0, \dots, l$ (when we denote $t_{0}=0$ and $t_{l+1}=1$) have the same triangulation up to isotopy on $\mathbb{S}^2$. We denote $T(s)$ for $s \in (t_{i},t_{i+1})$ by $T_{i+1}$. By the definition of the generic braid, at each critical time \( t_i \), the Delaunay triangulation of the configuration undergoes a flip, that is, $T_{i+1}$ is obtained from $T_{i}$ by a flip. Therefore, one can obtain a finite sequence of triangulations
\[
T_1 \rightarrow T_2 \rightarrow \cdots \rightarrow T_{\ell+1},
\] 
where $T_{i+1}$ is obtained from $T_{i}$ by a flip for $i=1,\dots, l$ and we call the sequence \textbf{a flip sequence} induced by the braid \( b \).



\begin{definition}

    Let $T:=T_{1} \rightarrow T_{2} \rightarrow \cdots \rightarrow T_{l}$ be a sequence of triangulations and flips. Assume that $T_{i}\rightarrow T_{i+1}\rightarrow \cdots \rightarrow T_{j}$ can be a part of involution, far-commutativity and pentagon relation, that is, $T_{i} \rightarrow T_{i+1}\rightarrow \cdots T_{j} \rightarrow T_{a} \rightarrow T_{a+1} \rightarrow \cdots T_{b} \rightarrow T_{i}$ forms one of involution, far-commutativity and pentagon relation. Let $T'$ be the sequence of triangulations and flips obtained from $T$ by replacing
    $$T_{i} \rightarrow T_{i+1}\rightarrow \cdots \rightarrow T_{j}$$
    by
    $$T_{i}\rightarrow T_{b} \cdots\rightarrow T_{a+1} \rightarrow T_{a} \rightarrow T_{j}.$$
    We say that {\bf $T'$ is obtained from $T$ by applying involution, far-commutativity and pentagon relation for flips on triangulation.}

    Flips related by involution, far-commutativity, and pentagon moves correspond to isotopies of the braid. Two isotopic braids induce flip sequences differing only by such relation.
\end{definition}

\begin{st}\label{prop:braid-flips}
        Let $\beta$ and $\beta'$ be two $n$-strand braids. Let $$T:=T_{1} \rightarrow T_{2} \rightarrow \cdots \rightarrow T_{l}$$
        and
        $$T'=T'_{1} \rightarrow T'_{2} \rightarrow \cdots \rightarrow T'_{l'}$$
        be sequences of flips corresponding to $\beta$ and $\beta'$.
        If $\beta$ and $\beta'$ are isotopic, then $T'$ is obtained from $T$ by applying involution, far-commutativity and pentagon relation for flips on triangulations.
\end{st}

\begin{proof}
Let \( \beta \) and \( \beta' \) be two isotopic \( n \)-strand braids. By definition, there is a continuous family \( \{\beta_s\}_{s \in [0,1]} \) such that:
\begin{itemize}
    \item \( \beta_0 = \beta \), \( \beta_1 = \beta' \),
    \item For all but finitely many parameters \( s_1, \ldots, s_u \in [0,1] \), the braid \( \beta_s \) is \emph{generic} in the sense that for all but finitely many times \( t \in [0,1] \), the configuration of points \( \{x_1^s(t), \ldots, x_n^s(t)\} \subset S^2 \) is in general position.
    \item
     At each exceptional parameter value \( s_i \), the isotopy passes through a codimension-two degeneracy; that is, one of the following occurs:
\end{itemize}

        \begin{enumerate}
             \item Four points lie on a circle at a single moment in time \( t \),
             \item Two independent cocircular quadruples occur simultaneously,
            \item Five points lie on a common circle.
        \end{enumerate}

Without loss of generality, we may assume that for each $s_{i}$ only one codimension two degeneracy appears in $\beta_{s}$.
Let \( T^{s} := T^{s}_1 \to T_2 \to \cdots \to T^{s}_l \) be the sequence of triangulations corresponding to $\beta_{s}$ for $s \in [0,1] \backslash \{s_{1}, \cdots, s_{u}\}$. Note that codimension two degeneracies give rise to involution, far-commutativity and pentagon relation. It follows that for each $s_{i} \in \{s_{1}, \cdots, s_{u}\}$, $T^{s_{i} +\epsilon}$ is obtained from $T^{s_{i} -\epsilon}$ by one of involution, far-commutativity and pentagon relation. It is clear that $T_{u}$ and $T_{v}$, $u,v \in (s_{i}, s_{i+1})$ are same up to isotopy on the sphere and hence the proof is completed.

\end{proof}

\section{Main Theorem: Construction of the invariant}

Let $(X, \oplus, \otimes)$ be a tropical semifield.
Let \(n \geq 5\) be an integer and \(\beta\) be a braid of \(n\) strands. Let \(T := T_1 \to T_2 \to \cdots \to T_l\) be a sequence of flips on Delaunay triangulations of the sphere \(\mathbb{S}^2\) with \(n\) points corresponding to the braid \(\beta\). Each triangulation \(T_i\) in the sequence has \(N = 3n - 6\) edges. We assume that the edges of each triangulation are labeled with elements from the tropical semifield \(X\).

Let \(A_1 = (a_1, a_2, \dots, a_N) \in X^{\times N}\) be a fixed labeling of the edges of the triangulation \(T_1\), where \(a_i\) is the tropical label of the \(i\)-th edge of \(T_1\). For each subsequent triangulation \(T_{i+1}\), the tropical labeling \(A_{i+1} = (a_1', a_2', \dots, a_N') \in X^{\times N}\) is obtained from \(A_i = (a_1, a_2, \dots, a_N)\) by a flip of labels.

Thus, the sequence of labels \(A_1, A_2, \dots, A_l\) evolves as a result of repeated flips.

\begin{definition}
    The map
\[
f_n(\beta) := A_\ell \in X^{\times N}
\]
is defined to be the final tropical labeling after applying all flips in the sequence corresponding to the braid $\beta$.
\end{definition}
    
\begin{theorem}
The map $f_n : SPB_n \to X^{\times N}$ is a braid invariant up to flip equivalence. That is, if two braids $\beta$ and $\beta'$ are isotopic, then $f_{n}(\beta)= f_{n}(\beta')$ in $X^{\times N}$.
\end{theorem}

\begin{proof}
    Let $\beta$ and $\beta'$ be isotopic braids. Let $$T:=T_{1} \rightarrow T_{2} \rightarrow \cdots \rightarrow T_{l}$$
        and
        $$T'=T'_{1} \rightarrow T'_{2} \rightarrow \cdots \rightarrow T'_{l'}$$
        be sequences of flips corresponding to $\beta$ and $\beta'$. By Proposition~\ref{prop:braid-flips} $T$ is obtained from $T'$ by applying involution, far-commutativity and pentagon relation for flips on triangulation. In Section \ref{subsec:consistency-labels} we showed that flips on labels are cosistent with involution, far-commutativity and pentagon relation for flips on triangulation. Since $f_{n}(\beta)$ and $f_{n}(\beta')$ are labelings of $T_{l}$ and $T'_{l'}$ obtained from $T$ and $T'$, it follows that $f_{n}(\beta) = f_{n}(\beta')$ in $X^{\times N}$. \qedsymbol
\end{proof}

\section{Acknowledgments}
I would like to express my sincere gratitude to my scientific supervisor, Manturov Vassily Olegovich, for his invaluable guidance, insightful advice, and continuous support throughout the preparation of this work. His expertise and encouragement have been instrumental in shaping the direction and quality of this paper. I am also deeply thankful to Kim Seongjeong for his generous assistance and helpful discussions, which significantly contributed to my understanding and his assistance in the preparation of this paper.

\end{document}